# On the convection-dispersion equation for a finite domain: Third-type boundaries as a necessary condition of the conservation law

W. J. Golz

*Department of Civil & Environmental Engineering, Louisiana State University, Baton Rouge, Louisiana 70803*

**Abstract**

This paper resolves a longstanding discussion of a mathematical problem important in contaminant hydrogeology and chemical-reaction engineering, by discussing the foundations for a conceptual model of a dilute miscible solute undergoing longitudinal convection and dispersion with moderate rates of appearance and disappearance in a finite continuum. It is demonstrated that: (i) Hulburts conditions (a first-type entrance with a third-type exit) fail to satisfy overall mass conservation; (ii) the conditions of Wehner and Wilhelm which reduce to those of Danckwerts (a third-type entrance with a zero-gradient exit) satisfy overall mass conservation yet fail to satisfy internal consistency with the governing equation; (iii) only third-type boundaries simultaneously satisfy internal consistency and overall mass conservation which are, respectively, a necessary and sufficient condition for any solution to the governing equation. This result is extensible to quite general governing equations since the boundary conditions are shown to be independent of the fate mechanisms.

*Keywords:* Biotransport; Ecohydrology; Environmental transport; Hydrodynamics; Insitu sediment caps; Soil columns

## 1. Introduction

Investigations involving porous media present a unique challenge because solute behavior can change dramatically at a boundary, and in many problems of practical interest it is impossible to measure interior concentrations without disturbing natural flow lines and thereby introducing error. As a consequence, empirical observations typically furnish only influent and effluent concentrations. Thus, our understanding of solute behavior in systems such as hydrogeologic strata or packed chemical reactors is often wholly dependent upon idealizations codified in conceptual models, and these models are quite sensitive to their description of mass transfer across the boundaries.

The diffusion of a miscible solute has long been recognized as an important phenomenon, and the earliest authors to approach general theoretical descriptions of this problem did so using common practice and intuition as their guide to boundary conditions. In 1944, Hulburt summarily assumed a first-type entry and postulated that, because typical reactor design would call for a reaction to reach its completion upstream of the outflow boundary, a zero gradient was the proper exit condition [6]. In 1953, Danckwerts explicitly satisfied the mass balance for a finite column by adopting a third-type entrance [3]. After noting that mass conservation led also to a third-type exit, Danckwerts went on to state that a positive exit gradient would require the reactant to pass through an interior minimum, while a negative gradient would lead to a condition where the concentration within the outlet was smaller than that in the external effluent stream.[1] Danckwerts regarded those causes of non-zero gradients as intuitively objectionable and sufficient for the dismissal of the third-type exit in favor of a zero gradient which, we note, implied that resident concentrations remained continuous across the lower boundary. In 1956, Wehner and Wilhelm extended the description of a finite reactor to a vessel having dispersive fore and aft reservoirs [10]. After equating third-type conditions at the upper and lower reactor-reservoir boundaries, Wehner and Wilhelm obtained a sufficient number of auxiliary conditions by introducing the supplemental assumption that resident concentrations remained continuous across the boundaries which, coincidentally, yielded a solution whose reactor concentrations were identical to those of Danckwerts.

Although the continuity of a resident concentration across a boundary has been a central question in an extensive interdisciplinary dialogue, the key qualities of that discussion can be illustrated by following the presentation of opposing postulates though several representative papers: In 1962, based upon the implied assumption that a Wehner-Wilhelm solution was correct, Bischoff and Levenspiel concluded that a solution employing a third-type entry in a semi-infinite domain incorrectly predicted the Peclet number for a finite column [1].[2] In 1978, Kreft and Zuber considered solutions for a semi-infinite domain and showed that a third-type entrance yielded the proper description for a resident concentration [7, table 1, eq. (6)], from which a proper Peclet number immediately follows. In that same paper, Kreft and Zuber went on to show that a first-type entrance described a flux concentration [7, table 1, eq. (7)] and that an expression for the resident concentration could be transformed into one describing flux [7, table 3, code T1].

---

[1] For Danckwerts' special case, where there is no production within the reactor, the dismissal of a positive exit gradient is correct, but the conservation law requires a negative gradient for any non-zero concentration.

The facts established by Kreft and Zuber were used by van Genuchten and Parker in 1984 as a basis for evaluating the effects of boundary conditions on short laboratory soil cores [9]. Van Genuchten and Parker showed that a Danckwerts solution subjected to a flux transformation failed to coincide with the curve for a Hulburts solution [9, fig. 2], lending credence to their contention that a zero gradient was a convenient assumption. Van Genuchten and Parker went on to demonstrate that a Danckwerts solution satisfied the same overall mass-balance as a semi-infinite solution for a first-type entrance [9, eq. [22]], yet the effluent curves for these two solutions were later shown to differ at early times when the Peclet number was small [9, fig. 3]. The Danckwerts solution should have furnished an expression with satisfactory convergence over the range of times and Peclet numbers plotted (see [9, fig. 3] then cf. [2, fig. 1]), as should the semi-infinite solution. This inconsistency illustrates the difficulty of drawing conclusions from contrapositives based upon graphical comparisons, absent a critical numerical analysis.

In a 1992 paper, Parlange et al. provided an extended argument in favor of regarding the resident concentration within the exit as lying somewhere between the discontinuity described by a semi-infinite solution with a third-type entry and the complete continuity described by a Danckwerts exit [8]. Their arguments, however, rest upon the dismissal of an infinitesimal element, which they accept at the entry, in favor of an outflow boundary layer which admits preferential flow paths and separable effects from the individual components of hydrodynamic dispersion. As we shall see, admitting those types of inhomogeneities contravenes the continuum assumptions required to write the governing CDE.

In 2001, Golz and Dorroh demonstrated that a Danckwerts exit was not generally valid because it implied a zero exit concentration [5]. That result is generalized in this paper where we use the conservation law as the basis for a fundamental theorem to show that third-type boundaries are a necessary and sufficient condition for any solution to the stated CDE (i.e., If a solution does not satisfy third-type boundaries, then it cannot be a solution to the CDE).[3]

---

[3] This theorem and its corollary were proved, albeit with slightly less generality, in 2003 by Golz whom furnished a particular solution to the governing equation and then demonstrated that self consistency required that solution to satisfy third-type boundaries [see 4, sec. 7.2 and app. D].

---

[2] It is customary to define the Peclet number for this type of problem as the product of fluid velocity and length of flow domain divided by the dispersion constant, i.e., $v\ell/D$.